# Optimization of On-condition Thresholds for a System of Degrading Components with Competing Dependent Failure Processes

Sanling Song, Nooshin Yousefi, David W. Coit
*Department of Industrial & Systems Engineering, Rutgers University, Piscataway, NJ, 08854*
Qianmei Feng
*Department of Industrial Engineering, University of Houston, Houston, TX 77204*

**Abstract** - An optimization model has been formulated and solved to determine on-condition failure thresholds and inspection intervals for multi-component systems with each component experiencing multiple failure processes due to simultaneous exposure to degradation and shock loads. In this new model, we consider on-condition maintenance optimization for systems of degrading components, which offers cost benefits over time-based preventive maintenance or replace-on-failure policies. For systems of degrading components, this can be a particularly difficult problem because of the dependent degradation and dependent failure times. In previous research, preventive maintenance and periodic inspection models have been considered; however, for systems whose costs due to failure are high, it is prudent to avoid the event of failure, i.e., we should repair or replace the components or system before the failure happens. The determination of optimal on-condition thresholds for all components is effective to avoid failure and to minimize cost. Low on-condition thresholds can be inefficient because they waste components life, and high on-condition thresholds are risky because the components are prone to costly failure. In this paper, we formulated and solved a new optimization model to determine optimal on-condition thresholds and inspection intervals. In our model, when the system is inspected, all components are inspected at that time. An inspection interval may be optimal for one component, but might be undesirable for another component, so the optimization requires a compromise. The on-condition maintenance optimization model is demonstrated on several examples.

*Keywords*: Multiple dependent competing failure processes, Degradation, Gamma process,

On-condition thresholds

## 1. Introduction

An effective maintenance policy maintains the system by achieving high safety and low cost, both of which are critical issues concerned in many modern industries [1]. Due to the inevitable deterioration of many components and system, systems may fail. To restore a failed system is often time-consuming and costly. Periodic and frequent inspection and repair/replacement can reduce the probability of deterioration and failure; however, it also incurs potentially excessive maintenance cost [2]. High quality operational performance and low maintenance cost can then become two conflicting objectives.

There has been much noteworthy research on reliability analysis for system subject to dependent failure processes, and accordingly, different maintenance policies have been considered [3-8]. For systems whose penalty cost due to downtime is high, detecting the component status and assisting in repair/replacement decision-making before system failure, leads to low risk of failure, and subsequently, lower maintenance cost. There have been previous studies on developing periodic inspection models for degrading system with components sharing dependent degradation and dependent failure time [9-11]; however, those maintenance models are generally not combined with on-condition thresholds for components. In this paper, we develop an optimization model to determine on-condition thresholds and inspection intervals for multi-component systems with each component experiencing multiple failure processes.

We initially present a reliability model for systems in which failure processes for each component are dependent and failure times for all components are dependent [6]. Second, we introduce a working principle for defining the on-condition threshold and system status. A periodic inspection maintenance policy is selected so that the decision-making depends on the on-condition thresholds for all components. Finally, a maintenance cost rate model is developed and minimized. In this model, system inspection interval and component on-condition thresholds are the decision variables. The new model offers cost benefits and performance improvement over time-based preventive maintenance or replace-on-failure policies.

The paper is organized as follows. Section 2 introduces and summarizes the relevant research previously done on reliability and maintenance policies for multi-component systems with multiple failure processes and provides the details of two failure processes. After presenting the system reliability model, Section 3 introduces on-condition threshold and defines system status related to the on-condition thresholds. Section 4 describes the maintenance policy and cost rate optimization model based on system inspection interval and component on-condition thresholds. System examples are shown in Section 5 to illustrate the developed reliability and maintenance models.

The notation used in formulating the reliability and maintenance models is listed as follows:

$N(t)$ = number of shock loads that have arrived by time $t$;
$n$ = number of components in a series or parallel system;
$\lambda$ = arrival rate of random shocks;
$D_i$ = threshold for catastrophic/hard failure of $i^{th}$ component;
$W_{ij}$ = size/magnitude of the $j^{th}$ shock load on the $i^{th}$ component;
$F_{Wi}(w)$ = cumulative distribution function (cdf) of $W_i$;
$H_i^1$ = critical wear degradation failure threshold of the $i^{th}$ component (a fixed parameter);
$H_i^2$ = on-condition threshold of the $i^{th}$ component (a decision variable);
$X_i(t)$ = wear volume of the $i^{th}$ component due to continuous degradation at $t$;
$X_{Si}(t)$ = total wear volume of the $i^{th}$ component at $t$ due to both continual wear and instantaneous damage from shocks;
$Y_{ij}$ = damage size contributing to soft failure of the $i^{th}$ component caused by the $j^{th}$ shock load;

| | | |
|---|---|---|
| $S_i(t)$ | = | cumulative shock damage size of the $i^{th}$ component at $t$; |
| $\alpha_i(t), \beta_i$ | = | shape and scale parameter for gamma degradation process for component $i$; |
| $G_i(x_i,t)$ | = | cumulative distribution function (cdf) of $X_i(t)$; |
| $F_{X_i}(x_i,t)$ | = | cdf of $X_{S_i}(t)$; |
| $f_{Y_i}(y)$ | = | probability density function (pdf) of $Y_i$; |
| $f_{Y_i}^{<k>}(y)$ | = | pdf of the sum of $k$ independent and identically distributed (i.i.d.) $Y_i$ variables |
| $f_T(t), F_T(t)$ | = | pdf and cdf of the failure time, $T$; |
| $F_T^{H^1}(t)$ | = | cdf of the failure time $T$ for the whole system considering critical failure threshold; |
| $F_T^{H^2}(t)$ | = | cdf of the time when an on-condition threshold is reached; |
| $C(t)$ | = | cumulative maintenance cost by time $t$; |
| $E[TC]$ | = | expected value of the total maintenance cost of the renewal cycle, $TC$; |
| $\tau$ | = | periodic inspection interval; |
| $CR(\tau)$ | = | average long-run maintenance cost rate of the maintenance policy; |
| $E[K]$ | = | expected renewal cycle length, $K$ of the maintenance policy; |
| $E[N_I]$ | = | expected number of inspections $N_I$; |
| $E[\rho]$ | = | expected system downtime (the expected time from a system failure to the next inspection when the failure is detected); |
| $C_R$ | = | replacement cost per unit; |
| $C_I$ | = | cost associated with each inspection; |
| $C_\rho$ | = | penalty cost rate during downtime per unit of time; |

## 2. Component and System Reliability Based on Degradation Analysis

Significant and meaningful prior research has been done on reliability and maintenance policies for systems with degradation, shocks, and independent or dependent failure processes. In this new system model, we extend previously developed models and research to develop a new maintenance optimization model to determine optimal component on-condition thresholds and system inspection interval.

### 2.1 Research on Reliability and Maintenance Models

There is a significant literature already dedicated to reliability analysis for system subject to multiple failure processes. Song et al. [6] studied the reliability of multi-component systems with each component experiencing multiple failure processes. This work was an extension of Peng et al. [35]. Chatwattanasiri et al. [12] then proposed a reliability model for a system of components with multiple competing and dependent failure processes when the future conditions are uncertain. Jiang et al [13] further studied reliability of systems subjected to multiple competing dependent failure processes with changing dependent failure thresholds.

There have been other studies for systems experiencing degradation processes and external random shocks. Wang and Pham [14] developed a model considering the dependent relationship between random shocks and degradation processes by a time-scaled covariate factor. Rafiee et al. [15] studied reliability for systems subject to dependent competing failure processes with a changing degradation rate according to particular random shock patterns. Jiang et al. [16] developed reliability model for system experiencing stochastic degradation process and a random shock process, with shock effects falling into distinct zones.

Different maintenance policies for degrading systems with a single component or multiple components have also been extensively studied in the literature [17]. Bian and Gebraeel [18] proposed a stochastic model for the degradation processes of components and estimated residual lifetime distribution of each component. Levitin and Lisnianski [19] studied a preventive maintenance optimization problem for multi-state systems, which have a range of performance levels. Tsai [20] proposed a preventive maintenance model for systems with deteriorating components. A simple preventive maintenance task is to restore the degraded component to some level of the original condition and a preventive replacement task is to replace the aged component by a new one or to restore to an as-new state.

Li and Pham [21] developed a generalized condition-based maintenance model subject to multiple competing failure processes including two degradation processes and random shocks, in which the preventive maintenance thresholds for degradation processes and inspection sequences are the decision variables. Grall et al. [22] focused on the analytical modeling of a condition-based inspection/replacement policy for a stochastically and continuously deteriorating single-unit system, in which both the replacement threshold and the inspection schedule are considered as decision variables for this maintenance problem. Tian and Liao [23] investigated condition-based maintenance policies of multi-component systems based on a proportional hazards model, where economic dependency exists among different components subject to condition monitoring. Perez et al. [24] proposed a method for scheduling the maintenance in a wind farm with multiple turbines each has multi components.

Jardine et al. [25] has performed recent research in diagnostics of mechanical systems implementing condition-based maintenance with emphasis on models, algorithms and technologies for data processing and maintenance decision-making. Optimizing condition-based maintenance for equipment subject to vibration has been studied by Jardine et al. [26]. Zhu et al. [27] considered a maintenance model for systems with degradation which are continuously monitored, and units are immediately repaired when failure happened. This process was repeated until a predetermined time was reached for preventive maintenance to be performed. Wang and Pham [28] studied a multi-objective maintenance optimization problem for systems subject to dependent competing risks of degradation wear and random shocks. The number of preventive maintenance actions until replacement and the initial preventive maintenance interval were determined by simultaneously maximizing the asymptotic system availability, and minimizing the system cost rate using the fast elitist Non-dominated Sorting Genetic Algorithm (NSGA).

Ko and Byon [29] used asymptotic theory to analytically solve the large scale maintenance optimization problem when the maintenance set up cost is higher than repair cost. Abdul-Malak and Kharoufeh [30] developed a Markov decision process model to find the optimal replacement strategy for a system of multiple components in a shared environment. Wang et al. [31] considered a multi-phase inspection schedule for a system that its degradation process is divided into more than two stages. Interaction between failure rates of units are considered for a two-unit system which is subjected to external shocks by Sung et al.[32].

**2.2 Review on Gamma Process**

In this paper, it is considered that each component degrades so that irreversible damage gradually occurs, and the degradation model is monotonically increasing. In this case, it is appropriate to use the gamma process to model the degradation path. A thorough review of the gamma process model and its applications can be found in Van Noortwijk [33]. For our applications, the gamma process with positive shape parameter is linear in $t$, with shape parameter $\alpha(t) = \alpha t$ and scale parameter $\beta$ is a continuous time stochastic process with the following properties:

- It starts from 0 at time 0, i.e., $X(0) = 0$
- $X(t)$ has independent increment
- for $t > 0$ and $s > 0$, $X(t) - X(s) \sim gamma(\alpha(t-s), \beta)$.

In fact, the probability density function of degradation process for each component $X_i(t) - X_i(s)$ is given by:

$$g(x; \alpha_i(t-s), \beta_i) = \frac{\beta_i^{\alpha_i(t-s)} x^{\alpha_i(t-s)-1} \exp(-\beta_i x)}{\Gamma(\alpha_i(t-s))} \tag{1}$$

where $\alpha_i(t)$ and $\beta_i$ are the shape parameter and scale parameter for component $i$.

Caballé et al. [34] proposed a condition-based maintenance strategy for a system that its degradation process follows a nonhomogenous Poisson process and its growth is modeled by gamma process.

## 2.3 Component Reliability with Competing Dependent Failure Processes

In this paper, we consider systems where each component can fail due to two competing dependent failure processes that share the same shock process; soft failure process and hard failure process [1, 2], as depicted in Figure 1. Each component in the system degrades with time, and when a shock arrives, if damage is greater than hard failure threshold, catastrophic failure occurs. If the component survives the shock, total degradation containing both pure degradation and additional incremental degradation caused by shock damage is greater than a defined soft failure threshold level, soft failure occurs. The two failure processes are competing and dependent.

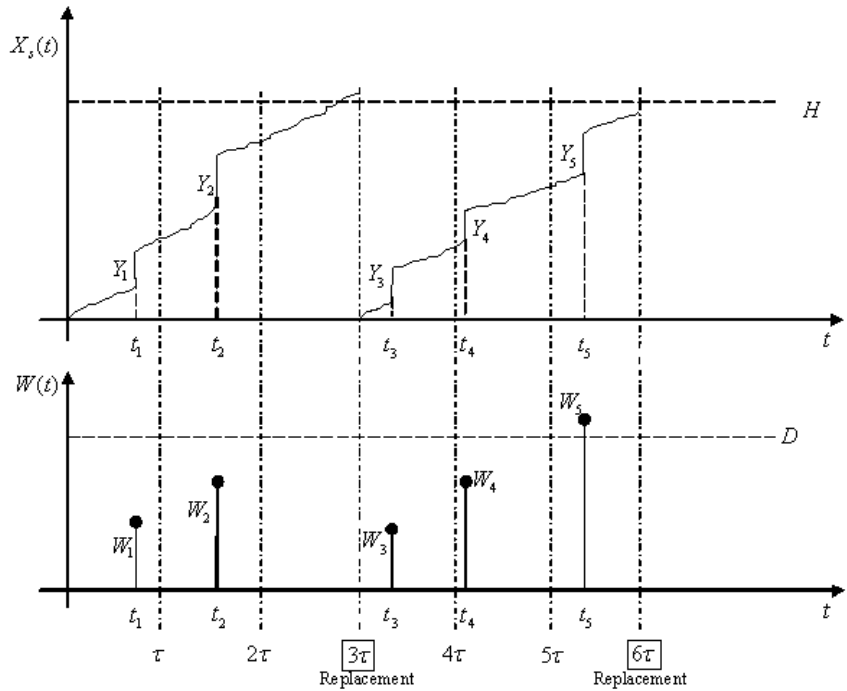

Figure 1. Two dependent and competing failure processes for a component

(a) soft failure process and (b) hard failure process [2]

Specific assumptions used for the reliability and maintenance modeling in this paper are as follows [1, 2]:

1. Soft failure occurs for the $i$th component when the total degradation of that component

exceeds its critical threshold level $H_i^1$. Component degradation is accumulated by both continuous degradation over time and cumulative incremental damage due to random shocks.

2. When the shock size exceeds the hard failure threshold of the components $D_i$, hard failure occurs of that component.
3. Random shocks arrive as a Poisson process.
4. The model is for systems that are packaged and sealed together, making it impossible or impractical to repair or replace individual components within the system, e.g., MEMS.
5. For the maintenance policy, the system is inspected at periodic intervals and no continuous monitoring is performed. Replacements are assumed to be instantaneous and perfect.
6. Upon an inspection for a series system, if the overall degradation of all $n$ components is lower than $H_i^2$, the system is within the high safety level area, and nothing is done.
7. If the degradation of any component is between $H_i^1$ and $H_i^2$, the system does not fail, but we consider it prone to high failure risk, and the system is replaced with a new one preventively.
8. If the system fails, that is, the total degradation of any component is higher than $H_i^1$ before the specified inspection interval, it is not immediately detected and not replaced until the next inspection. There is penalty cost per time associated with the failure of system during downtime, e.g., cost associated with loss of production, opportunity costs, etc.

In this paper, we develop an optimization model to determine on-condition failure thresholds and inspection intervals for complex multi-component system with each component experiencing multiple failure processes due to simultaneous exposure to degradation and shock loads. Two failure processes for each component are dependent, and failure times for all components are also dependent. Component hard failures occur when a shock load exceeds thresholds. Figure 1(b) shows that component $i$ may fail when damage from a shock exceeds $D_i$. $W_{ij}$ is the shock size and it is an *i.i.d.* random variable with some defined distribution which is assume in this paper as a normal distribution, $W_{ij} \sim N(\mu_{W_i}, \sigma_{W_i}^2)$, although this is not a restriction for our model. We can obtain the probability that the $i^{\text{th}}$ component survives a shock [2]:

$$P_{Li} = P(W_{ij} < D_i) = F_{Wi}(D_i) = \Phi\left(\frac{D_i - \mu_{W_i}}{\sigma_{W_i}}\right) \text{ for } i = 1, 2, \ldots, n, \tag{2}$$

where $\Phi(\cdot)$ is the cdf of a standard normal random variable.

As shown in Figure 1(a), total degradation of the $i^{\text{th}}$ component can be accumulated as $X_{S_i}(t) = X_i(t) + S_i(t)$, and when $X_{S_i}(t) > H_i^1$, soft failure occurs. Conditioning on the number of shocks and using a convolutional integral of $X_{S_i}(t)$, we can obtain the probability that component $i$ does not experience soft failure before time $t$ as follow:

$$P(X_{S_i}(t) < H_i^1) = F_{X_i}(H_i^1, t) = \sum_{m=0}^{\infty} P\left(X_i(t) + \sum_{j=1}^{m} Y_{ij} < H_i^1\right) \frac{\exp(-\lambda t)(\lambda t)^m}{m!} \tag{3}$$

$$P(X_{S_i}(t) < H_i^1) = F_{X_i}(H_i^1, t) = \sum_{m=0}^{\infty} \left(\int_0^{H_i^1} G_i(H_i^1 - u, t) f_{Y_i}^{<m>}(u) du\right) \frac{\exp(-\lambda t)(\lambda t)^m}{m!} \tag{4}$$

$X_i(t)$ follows a gamma process, so $G_i(\cdot)$ is the cdf for a gamma distribution. It is convenient for $Y_i$ to be gamma or normal distributed because the sum of $m$ iid gamma random variables is also gamma, and the sum of $m$ iid normal random variables is normal. In Song et al. [6], the assumption was made that $Y_i$ was normally distributed, while in this paper, we assume $Y_i$ is gamma distributed, but this is not a restriction.

## 2.4 Reliability Analysis for Multiple Components System with MDCFP

Our example system configuration is a series system, in which a component fails when either of the two dependent and competing failure modes occurs, and all components in the system behave similarly. Song at al. [6] developed a multi-component system reliability model when each component experiencing multiple failure processes due to each component degradation and external shock loads. The reliability of this series system can be obtained, since the system fails when the first component fails. The concepts described in this paper can be extended to other system configurations as well.

Figure 2 shows a series system with *n* components. The reliability of this series system at time *t* is the probability that each component survives each of the *N(t)* shock loads ($W_{ij} < D_i$ for $j = 1, 2, \ldots$) and the total degradation of each component is less than the soft failure threshold level ($X_{Si}(t) < H_i^1$ for all *i*).

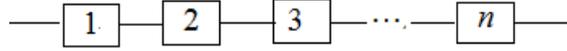

Figure 2. Series system example

In this model, shocks arriving at random time intervals are modeled as a Poisson process. When the system receives a shock (at rate $\lambda$), all components experience a shock. If we consider the component survival probabilities conditioned on the number of shocks, then the failure processes for all components become independent for a fixed number of shocks. System reliability function can be derived for the general case for a series system as follows [6]:

$$R(t) = \sum_{m=0}^{\infty} \prod_{i=1}^{n} \left[ P(W_i < D_i)^m P\left( X_i(t) + \sum_{j=1}^{m} Y_{ij} < H_i^1 \right) \right] \frac{\exp(-\lambda t)(\lambda t)^m}{m!} \qquad (5)$$

Using convolution integral, it can be obtained as follow [6]:

$$R(t) = \sum_{m=0}^{\infty} \prod_{i=1}^{n} \left[ P(W_i < D_i)^m \int_0^{H_i^1} G_i(H_i^1 - u, t) f_{Y_i}^{<m>}(u) du \right] \frac{\exp(-\lambda t)(\lambda t)^m}{m!} \qquad (6)$$

## 3. Operational Principle of the On-Condition Rule

For systems whose costs associated with failure are high, it is advantageous to repair or replace the components or system before the failure occurs. The concept of condition monitoring and on-condition thresholds for the components is used to evaluate and measure system status, and therefore, increase the opportunity to detect the components' critical and degraded situation and to avoid costly failure events. Maintenance optimization is based on reliability modeling of system subject to dependent and competing failure processes. The maintenance optimization is challenging because of the dependent degradation and dependent failure times among all components.

**3.1 Definition of system status related to on-condition threshold**

For some systems, the cost and consequence of failure are excessive compared to comparable preventive repair cost, replacement cost or other kinds of cost. Therefore, it is prudent to prevent failure from occurring and replace the equipment at the earliest convenience after it has sufficiently aged, rather than allowing to fail and possibly cause more severe consequences. For the multi-component system considered in this research, the components are packaged and/or sealed together and it is reasonable or necessary to replace the whole system before the critical degradation thresholds are reached. On-condition rules provide the capability to measure system status and detect failures before they occur so that preventive maintenance can be performed. Based on defined rules, the implementation of a lower degradation threshold can be useful to avoid failure by providing criteria to detect the degradation status of the components.

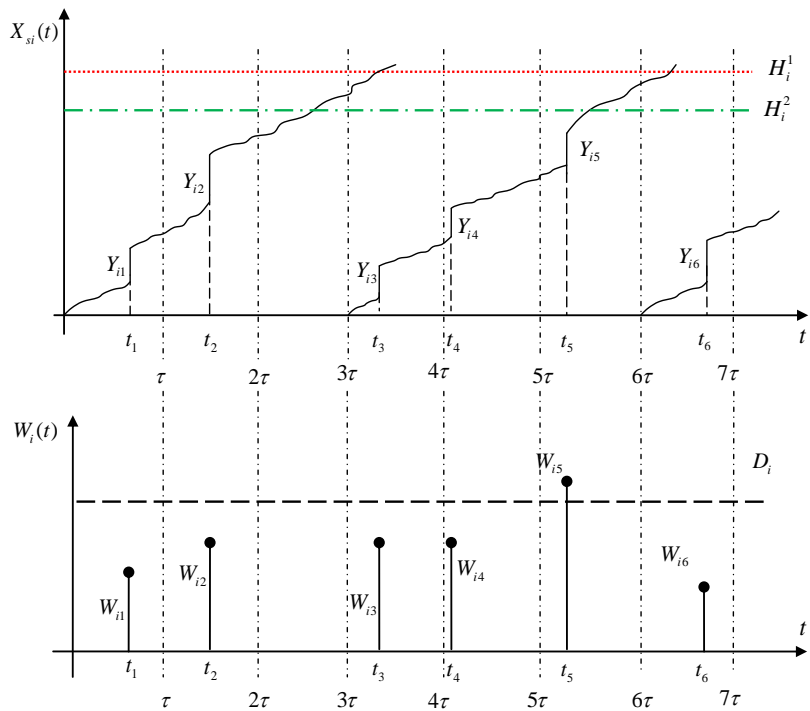

Figure 3. Two thresholds divide system status into three regions

$H_i^1$ is defined as the soft failure threshold for component *i* and $H_i^2$ is now defined as the on-condition threshold for component *i*, with $H_i^2 \leq H_i^1$. At each inspection time, we determine component condition for each component by inspection and compare it to a threshold. The action taken depends on a selection of condition-based operational status data and the defined maintenance condition rules. In Fig. 3, we can see given a fixed on-condition threshold $H_i^2$ for component *i* (lower bar and dash line in soft failure process). We adopt rules related to this on-condition degradation threshold to define the component degradation state.

At each inspection interval, if no hard failure occurs, and at the same time, total degradation of the $i^{th}$ component is less than $H_i^2$, the on-condition threshold for $i^{th}$ component, we then consider the component is in the safe region. The safe region is defined as the combination of soft failure process and hard failure process both below their respective thresholds and this status is defined as event A shown in Table 1. If no hard or soft failure occurs and total degradation is between $H_i^2$ and $H_i^1$ for any component *i*, this component has not failed; however, probabilistically it may fail within a short period. This status can be described by the combination of soft failure process area between $H_i^2$ and $H_i^1$, and hard failure process area below the hard failure threshold, which is defined as event B as presented in Table 1. If there has been a hard failure or the total degradation of any component *i* is greater than $H_i^1$ (higher dash line in soft failure process), the system has failed. The status can be defined as the union of the soft failure process area above the red dashed line, and hard failure process area above black dashed line, and this status is defined as event C.

Table1. Component status defined with two soft failure thresholds and hard failure threshold

| A | component is in safe region | $P(A) = \sum_{m=0}^{\infty} P(W_i < D_i)^m \left( \int_0^{H_i^2} G_i(H_i^2 - u, t) f_{Y_i}^{<m>}(u) du \right) \frac{\exp(-\lambda t)(\lambda t)^m}{m!}$ |
|---|---|---|
| B | component is working, but probabilistically fails soon | $P(B) = \sum_{m=0}^{\infty} P(W_i < D_i)^m \left( \int_0^{H_i^1} G_i(H_i^1 - u, t) f_{Y_i}^{<m>}(u) du - \int_0^{H_i^2} G_i(H_i^2 - u, t) f_{Y_i}^{<m>}(u) du \right)$ $\times \frac{\exp(-\lambda t)(\lambda t)^m}{m!}$ |
| C | component fails | $P(C) = \sum_{m=0}^{\infty} \left( 1 - P(W_i < D_i)^m \int_0^{H_i^1} G(H_i^1 - u, t) f_Y^{<m>}(u) du \right) \frac{\exp(-\lambda t)(\lambda t)^m}{m!}$ |

The probability that component total degradation less than $x$ by time $t$ as $\Psi(x)$:

$$\Psi(x) = \sum_{m=0}^{\infty} \left( \int_0^x G_i(x - u, t) f_{Y_i}^{<m>}(u) du \right) \frac{\exp(-\lambda t)(\lambda t)^m}{m!} \qquad (7)$$

Considering the safe region case for example, conditioning on $m$ shocks arriving to the system by time $t$ with probability $\frac{\exp(-\lambda t)(\lambda t)^m}{m!}$, the probability of no hard failure is $P(W_i < D_i)^m$, and the probability that total degradation is less than $H_i^2$ is $\Psi(H_i^2) = \int_0^{H_i^2} G_i(H_i^2 - u, t) f_{Y_i}^{<m>}(u) du$. Combining both soft failure process and hard failure process, the probability for event $A$: the component is in safe region is:

$$P(A) = \sum_{m=0}^{\infty} P(W_i < D_i)^m \left( \int_0^{H_i^2} G_i(H_i^2 - u, t) f_{Y_i}^{<m>}(u) du \right) \frac{\exp(-\lambda t)(\lambda t)^m}{m!}.$$

Similarly, for event $B$, component $i$ is still working, but it may probabilistically fail within the next inspection interval, the probability of no hard failure considering on m shocks is $P(W_i < D_i)^m$ and the probability that total degradation is between $H_i^1$ and $H_i^2$ is $\int_{H_i^2}^{H_i^1} G_i(H_i^1 - u, t) f_{Y_i}^{<m>}(u) du$. Combining both soft failure process and hard failure process; we can obtain the probability for event $B$. For event $C$, either soft failure or hard failure occurs, we can also determine probability, which equals to one minus the probability that neither of these two failure happens. The policy is summarized in Table 1.

Given this reliability model for systems with each component experiencing multiple failure processes due to simultaneous exposure to degradation and shock loads, we can define a maintenance cost optimization objective function. The system is inspected periodically, and the condition of each component is observed and compared to a threshold. Upon an inspection, we replace the system with a new one, when we observe that a hard failure has occurred or total degradation is greater than the on-condition threshold for any component $i$.

The expected number of inspections $N_I$, for a vector on-condition thresholds $\mathbf{H}^2 = (H_1^2, H_2^2, ..., H_n^2)$ is given by,

$$E(N_I) = \sum_{k=1}^{\infty} k(F_T^{\mathbf{H}^2}(k\tau) - F_T^{\mathbf{H}^2}((k-1)\tau)) \tag{8}$$

$$F_T^{\mathbf{H}^2}(t) = 1 - \left(\sum_{m=0}^{\infty} \prod_{i=1}^{n} \left[ P(W_i < D_i)^m \int_0^{H_i^2} G_i(H_i^2 - u, t) f_{Y_i}^{<m>}(u) du \right] \frac{\exp(-\lambda t)(\lambda t)^m}{m!} \right) \tag{9}$$

From Fig. 4, we can observe that system downtime is the time duration between the time a failure occurs and the next time an inspection is performed, and a failure detected. Conditioning on the event that there is a failure at time $t$ between the $(k-1)^{th}$ and $k^{th}$ inspection $[(k-1)\tau, k\tau]$ with probability $\left(F_T^{\mathbf{H}^2}(k\tau) - F_T^{\mathbf{H}^2}((k-1)\tau)\right)$, and defining the failure time as $t$, the system downtime is $k\tau - t$. The expected value of system downtime or the expected time from a system failure to the next inspection when the failure is detected, can then be determined as $\int_{(k-1)\tau}^{k\tau} (k\tau - t) dF_T^{\mathbf{H}^1}(t)$. Summing over the probability that failure can occur in any inspection interval, we can obtain expected system downtime as follows:

$$E[\rho] = \sum_{k=1}^{\infty} E[\rho | N_I = k] P(N_I = k) = \sum_{k=1}^{\infty} \left( \left(F_T^{\mathbf{H}^2}(k\tau) - F_T^{\mathbf{H}^2}((k-1)\tau)\right) \int_{(k-1)\tau}^{k\tau} (k\tau - t) dF_T^{\mathbf{H}^1}(t) \right) \tag{10}$$

$$F_T^{\mathbf{H}^1}(t) = 1 - \left( \sum_{m=0}^{\infty} \prod_{i=1}^{n} \left[ P(W_i < D_i)^m \int_0^{H_i^1} G_i(H_i^1 - u, t) f_{Y_i}^{<m>}(u) du \right] \frac{\exp(-\lambda t)(\lambda t)^m}{m!} \right) \quad (11)$$

The expected time between two replacements or expected cycle length is

$$E[K] = \sum_{k=1}^{\infty} E[K \mid N_I = k] P(N_I = k) = \sum_{k=1}^{\infty} k\tau (F_T^{\mathbf{H}^2}(k\tau) - F_T^{\mathbf{H}^2}((k-1)\tau)) \quad (12)$$

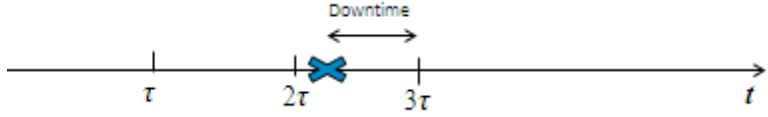

Figure 4. System downtime under periodic inspection maintenance policy

## 4. Condition-Based Maintenance Modeling and Optimization

With the on-condition rules stated in last section, a condition-based maintenance policy can be defined for the system with multiple components each exposed to two competing dependent failure processes. Condition-based maintenance offers the promise of enhancing the effectiveness of maintenance programs. We consider the practical case that when the penalty cost due to downtime is relatively higher than the corresponding preventive maintenance cost, so it is better to replace the whole system before the wear volumes of components reach their critical degradation thresholds. For some systems, it is best to just let them fail, but we are not considering those cases in this paper.

### 4.1 Description of Maintenance Model

On-condition degradation threshold can achieve our goal of replacing the system before failure by providing the criteria to detect the degradation of component beyond the on-condition threshold. If the on-condition threshold is too low and far away from the nominal threshold level, then we have to replace the whole system more frequently, and it results in extra cost due to the waste of system life. Alternatively, if the threshold is too high, then the system may fail before the next inspection leading to potentially expensive downtime cost. Therefore, on-condition degradation thresholds for all components and an inspection interval for the whole system are chosen to be decision variables in this maintenance optimization problem.

To evaluate the performance of the condition-based maintenance policy, we use an average long-run maintenance cost rate model, in which the periodic inspection interval $\tau$ for the whole system and on-condition thresholds $H_i^2$ for all components are the decision variables. At time $\tau$, and subsequent inspection intervals of time $\tau$, the entire assembled system is inspected. If the system is still operating satisfactorily with no component wear volume above the on-condition threshold, nothing is done. If degradation thresholds for all component are below the fixed critical degradation thresholds $H_i^1$ but some are above the on-condition threshold $H_i^2$, the whole system is replaced preventively. If there is a hard failure or at least one component's wear volume is above the critical degradation threshold $H_i^1$ prior to inspection, then the system is not replaced with a new one correctively until the next inspection. The average long-run maintenance cost per unit time can be evaluated by:

$$\lim_{t\to\infty}(C(t)/t) = \frac{\text{Expected maintenance cost between two replacements}}{\text{Expected time between two replacements}} = \frac{E[TC]}{E[K]} \quad (13)$$

where $TC$ is the total maintenance cost of a renewal cycle, and $K$ is the length of a cycle that takes a value of a multiple of $\tau$ [35]. The expected total maintenance cost is given as:

$$E[TC] = C_I E[N_I] + C_\rho E[\rho] + C_R \quad (14)$$

where $C_I$ is the cost of each inspection. $C_R$ is the replacement cost, $C_\rho$ is the penalty cost incurred during down time, and $\tau$ is the time interval for periodic inspection. Based on Eqs. (8) to (10), the average long-run maintenance cost rate is given as

$$CR(\tau, \mathbf{H}^2) = \frac{C_I \sum_{k=1}^{\infty} k\left(F_T^{\mathbf{H}^2}(k\tau) - F_T^{\mathbf{H}^2}((k-1)\tau)\right) + C_\rho \sum_{k=1}^{\infty} \left(\left(F_T^{\mathbf{H}^2}(k\tau) - F_T^{\mathbf{H}^2}((k-1)\tau)\right) \int_{(k-1)\tau}^{k\tau} (k\tau - t) dF_T^{\mathbf{H}^1}(t)\right) + C_R}{\sum_{k=1}^{\infty} k\tau (F_T^{\mathbf{H}^2}(k\tau) - F_T^{\mathbf{H}^2}((k-1)\tau))}$$

(15)

## 4.2 Maintenance Cost Optimization

For our maintenance optimization problem, if there are $n$ components in a series system, there are $n+1$ decision variables; namely $n$ on-condition thresholds for all components and the periodic inspection interval for the whole system. Our objective is to minimize maintenance cost rate, and constraints are that on-condition thresholds for all components should be less than or equal to their critical failure thresholds, and inspection interval should be a positive value. Therefore, our maintenance optimization problem can be formed as follows:

$$\begin{aligned}
\min \quad & CR(\tau, \mathbf{H}^2) \\
\text{s.t.} \quad & 0 \leq H_1^2 \leq H_1^1, \\
& 0 \leq H_2^2 \leq H_2^1, \\
& \quad ... \\
& 0 \leq H_n^2 \leq H_n^1, \\
& \tau \geq 0,
\end{aligned} \quad (16)$$

It is a difficult non-linear optimization problem but with continuous decision variables and a convex feasible region. For constrained nonlinear optimization problems, there are many available algorithms to obtain optimal solutions.

To solve the optimization problem, an interior point method is used (as implemented as the fmincon algorithm in the MATLAB optimization toolbox). The method consists of a self-concordant barrier function used to encode the convex set. It reaches an optimal solution by traversing the interior of the feasible region using one of two main types of steps at each iteration [36]. The algorithm first attempts to take a direct step within the feasible region to solve the Karush Kuhn Tucker (KKT) equations for the approximate problem by a linear approximation, which is also called a Newton step. By solving the KKT equations, we can get the direct step and the solution for the next iteration. If it cannot take a direct step, it attempts a conjugate gradient step, and minimizes a quadratic approximation to the approximate problem in a trust region, subject to linearized constraints. It does not take a direct step is when the problem is not locally convex near the current iteration. At each iteration, the algorithm decreases a merit function. We reach a new solution point after taking the step and start a new iteration. It continues until stopping criterion is met.

## 5. Numerical Examples

We consider two examples; the first one is a series system with four components where component 1 and 2 have the same parameters and component 3 and 4 are also the same. The parameters for reliability analysis are provided in Table 2. For this example, $Y_{ij}$ follows gamma distributions and $W_{ij}$ follows normal distributions. Without loss of generality, we assume that those parameters of component 1 and 2 are the same, and parameters of component 3 and 4 are the same. This is a conceptual example to demonstrate the reliability function and maintenance models. However, although the example is conceptual, $H_i^1$ and $D_i$ are estimated based on documented degradation trends [3]. In this part, we perform maintenance optimizations for both series system and all the individual components making up the system separately, and we discuss the results.

Table 2. Parameter values for multi-component system reliability analysis

| Parameter | component 1 & 2 | component 3 & 4 | Sources |
|---|---|---|---|
| $H_i^1$ | 0.00125 μm$^3$ | 0.00127 μm$^3$ | Tanner and Dugger [3] |
| $D_i$ | 1.5 Gpa | 1.4 Gpa | Tanner and Dugger [3] |
| $\alpha_i$ | 0.7 | 0.8 | Assumption |
| $\beta_i$ | 0.3 | 0.3 | Assumption |
| $\lambda$ | 2.5×10$^{-5}$ | 2.5×10$^{-5}$ | Assumption |
| $Y_{ij}$ | $Y_{ij} \sim gamma(\alpha_{Y_i}, \beta_{Y_i})$<br>$\alpha_{Y_i} = 0.4,\ \beta_{Y_i} = 1$ | $Y_{ij} \sim gamma(\alpha_{Y_i}, \beta_{Y_i})$<br>$\alpha_{Y_i} = 0.5,\ \beta_{Y_i} = 1$ | Assumption |
| $W_{ij}$ | $W_{ij} \sim N(\mu_{Wi}, \sigma_{Wi}^2)$<br>$\mu_{Wi}$ =1.2 GPa, $\sigma_{Wi}$ =0.2 GPa | $W_{ij} \sim N(\mu_{Wi}, \sigma_{Wi}^2)$<br>$\mu_{Wi}$ =1.22 GPa, $\sigma_{Wi}$ =0.18 GPa | Assumption |

First, we consider the maintenance policy for the whole series system with four components and a predetermined inspection interval, i.e., we inspect the whole system at one interval of $\tau$ and replace the system when the wear volume is above $H_i^2$ for any component. Choosing $C_I$=$1, $C_\rho$=$20000 and $C_R$=$100 and fixed interval of $\tau$=120 hours, we can find the minimum average long-run maintenance cost rate for system is $3.054×10^2$ and on-condition degradation threshold are $H_1^{2*}= H_2^{2*}$=0.0001556, $H_3^{2*}= H_4^{2*}$=0.0001370. By considering a shorter fixed inspection interval of $\tau$=24 hours, the minimum average long-run maintenance cost rate for system reduces to $2.2796×10^2$ and on-condition degradation threshold are $H_1^{2*}= H_2^{2*}$=0.0004637, $H_3^{2*}= H_4^{2*}$=0.0004204. When the system is inspected more frequently, we have higher on-condition degradation thresholds, i.e., closer to the failure threshold. Since the system status is detected more often, it can be replaced preventively, so on-condition degradation thresholds is closer to failure thresholds.

Inspection interval and on-condition threshold are two variables with an interesting trade-off. If we inspect component/system quite often, although we have higher on-condition thresholds, we may still have higher probability to detect system status and replace it preventively. If we inspect less often, we can compensate by defining a very low on-condition threshold value to achieving preventive maintenance.

The contribution of this paper is to now simultaneously determine the optimal on condition thresholds and inspection interval. The minimum average long run maintenance cost rate for the system is $1.9023×10^2$ found after 22 steps of iteration. The inspection interval is $\tau^*$=44.7129 hours, and on-condition degradation thresholds are $H_1^{2*}= H_2^{2*}$=0.0003055 and $H_3^{2*}= H_4^{2*}$=0.0002728. Figure 5 illustrates the iteration process of decision variables: inspection interval, on-condition degradation threshold for component 1 and 2, and on-condition degradation threshold for component 3 and 4. Figure 6 shows the iteration for our objective function, i.e., the system maintenance cost rate

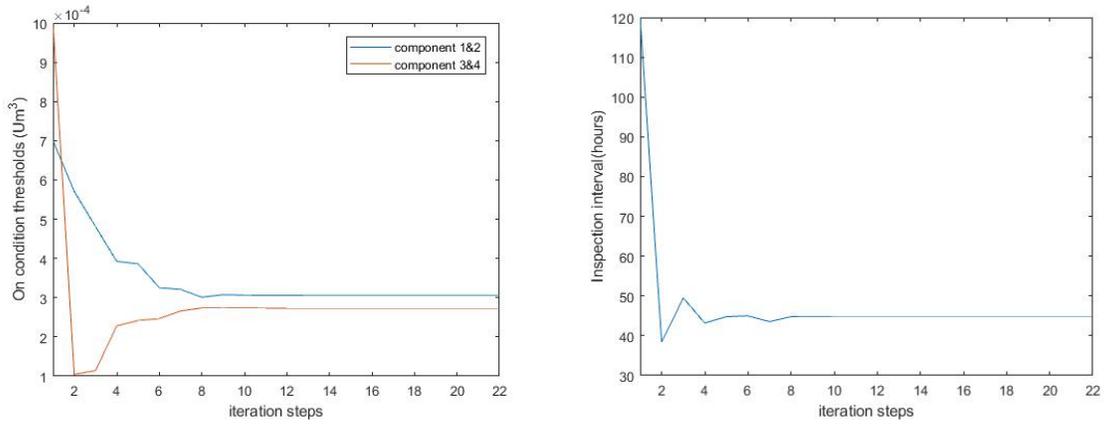

Figure 5. Iteration process for these decision variables: inspection interval $\tau^*$, and on-condition threshold for all components.

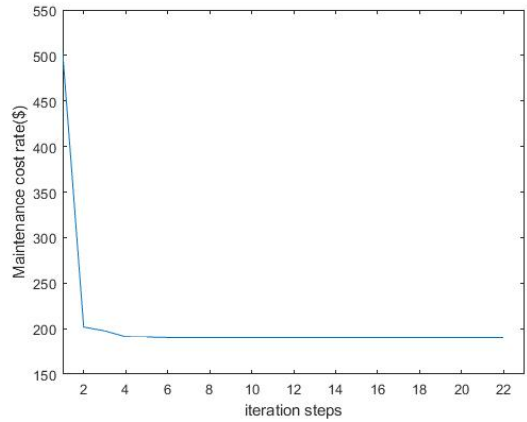

Figure 6. Iteration process of maintenance cost rate for system with four components

To further evaluate the results, we also consider an inspection and maintenance policy for the individual components. That is, we treat four components as individual systems, and inspect individual four components at their own inspection intervals. Since component 1 and 2 share the same parameter, the maintenance optimization for them are the same. We can find the minimum average long-run maintenance cost rate for component 1 and 2 as $\$1.367 \times 10^2$ after 20 steps of iteration, with a solution of the periodic inspection interval $\tau_{1,2}^*=65.044$ hours, and on-condition degradation threshold for components $H_1^{2*}= H_2^{2*}=0.0002465$. Figure 7 illustrates the iteration process of two decision variables, inspection interval and on-condition degradation threshold, for component 1 and 2. Figure 8 shows the iteration for our objective function, that is, the maintenance cost rate.

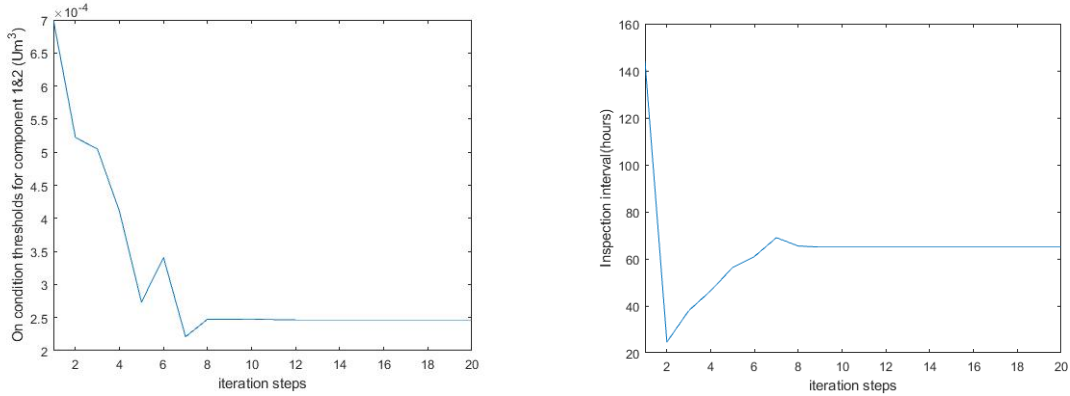

Figure 7. Iteration process two decision variables: inspection interval $\tau^*$, and on-condition threshold for component 1 and 2

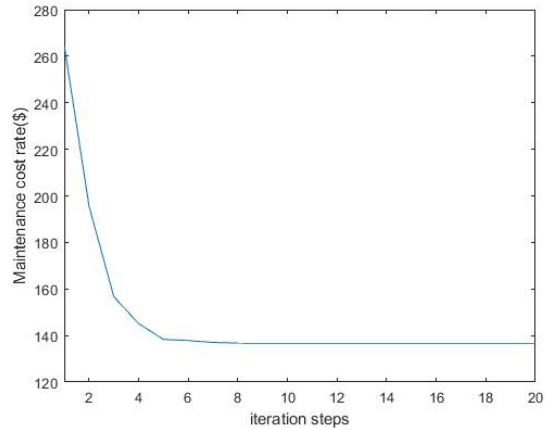

Figure 8. Iteration process of maintenance cost rate for component 1 and 2

Similarly, we inspect individual component 3 or component 4 at their own inspection intervals. The minimum average long-run maintenance cost rate for component 3 and 4 is $1.762×10^2$ after 13 steps of iteration, with the periodic inspection interval $\tau_{3,4}^*$=71.55 hours, and on-condition degradation threshold for components $H_3^{2*}= H_4^{2*}$=0.0002169. Figure 9 illustrates the iteration process of two decision variables: inspection interval and on-condition degradation threshold for component 3 and 4. Figure 10 shows the iteration for our objective function, that is, the maintenance cost rate.

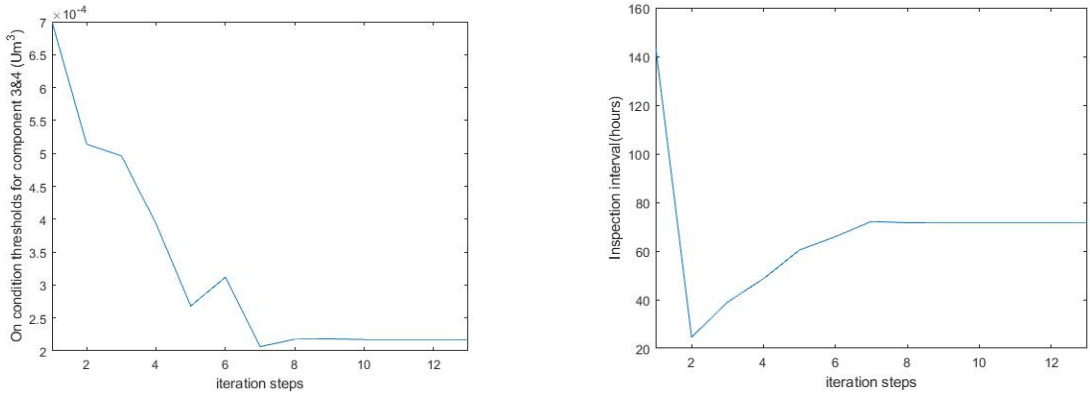

Figure 9. Iteration process two decision variables: inspection interval $\tau^*$, and on-condition threshold for component 3 and 4

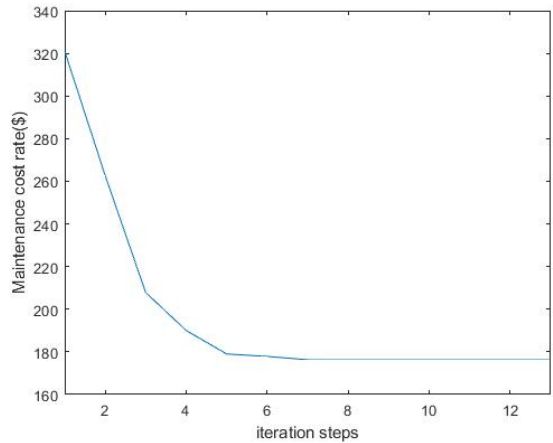

Figure 10. Iteration process of maintenance cost rate for component 3 and 4

We can observe that inspection intervals for either component 1 and 2 or component 3 and 4 are greater than the inspection interval for the series system, which means we have to compromise to inspect the system more frequently if we have more components in the system. Since time to failure for all components are different, and series system reliability is less than individual component reliability given any fixed time, we should inspect system more often to increase probability of avoiding failure and relative high downtime cost.

The second example is a series system with four different components. Table 3 presents the parameters of each component. Given the same cost $C_I$=$1, $C_\rho$=$20000 and $C_R$=$100, we find

the minimum average long-run maintenance cost rate for the system as $1.8356\times10^2$, which is obtained at periodic inspection interval $\tau^*$=49.86 hours, and on-condition degradation threshold for components are $H_1^{2*}$=0.0002904, $H_2^{2*}$=0.0002656, $H_3^{2*}$=0.0007362, $H_4^{2*}$=0.0012359. As the results illustrate, component 4 has the highest optimal on-condition threshold that is very close to its failure threshold. This is mainly because the degradation rate and shock load damage for component 4 is lower than other components which means its reliability is higher compared to all other three components. Accordingly, its optimal on-condition threshold is higher.

Table 3. Parameter values for multi-component system reliability analysis

| Parameter | component 1 | component 2 | component 3 | component 4 |
|---|---|---|---|---|
| $H_i^1$ | 0.00125 μm³ | 0.00127 μm³ | 0.0013 μm³ | 0.00128 μm³ |
| $D_i$ | 1.5 Gpa | 1.4 Gpa | 1.2 Gpa | 1.45 Gpa |
| $\alpha_i$ | 0.7 | 0.8 | 0.6 | 0.2 |
| $\beta_i$ | 0.3 | 0.3 | 0.25 | 0.25 |
| $\lambda$ | $2.5\times10^{-5}$ | | | |
| $Y_{ij}$ | $Y_{ij} \sim gamma(\alpha_{Y_i}, \beta_{Y_i})$ $\alpha_{Y_i} = 0.45, \beta_{Y_i} = 1$ | $Y_{ij} \sim gamma(\alpha_{Y_i}, \beta_{Y_i})$ $\alpha_{Y_i} = 0.5, \beta_{Y_i} = 1$ | $Y_{ij} \sim gamma(\alpha_{Y_i}, \beta_{Y_i})$ $\alpha_{Y_i} = 0.48, \beta_{Y_i} = 1$ | $Y_{ij} \sim gamma(\alpha_{Y_i}, \beta_{Y_i})$ $\alpha_{Y_i} = 0.4, \beta_{Y_i} = 1$ |
| $W_{ij}$ | $W_{ij} \sim N(\mu_{Wi}, \sigma_{Wi}^2)$ $\mu_{Wi}$=1.2 GPa, $\sigma_{Wi}$=0.22 GPa | $W_{ij} \sim N(\mu_{Wi}, \sigma_{Wi}^2)$ $\mu_{Wi}$=1.22 GPa, $\sigma_{Wi}$=0.18 GPa | $W_{ij} \sim N(\mu_{Wi}, \sigma_{Wi}^2)$ $\mu_{Wi}$=1.23 GPa, $\sigma_{Wi}$=0.15 GPa | $W_{ij} \sim N(\mu_{Wi}, \sigma_{Wi}^2)$ $\mu_{Wi}$=1.2 GPa, $\sigma_{Wi}$=0.2 GPa |

## 6 Conclusions

In this paper, we propose a maintenance optimization model to determine on-condition failure thresholds and inspection intervals for systems with dependent degradation and dependent component failure times. For systems whose penalty cost due to downtime is high, this on-condition maintenance policy offers cost benefits over time-based preventive maintenance or replace-on-failure, because on-condition threshold increases the liklihood to detect system critical status and prevent failures. In this maintenance policy, the periodic inspection interval for the whole system and on-condition thresholds for all components are decision variables, and system maintenance cost rate is our optimization objective. The average long-run maintenance cost rate is evaluated and optimized. Interior point algorithm in MATLAB toolbox fmincon is used to solve the optimization problem. Numerical examples are provided and the results are discussed.

**Acknowledgments**

This study was based in part upon work supported by USA National Science Foundation (NSF) grants CMMI-0970140 and CMMI-0969423.

systems using proportional hazards model," *Reliability Engineering & System Safety*, vol. 96, no. 5, pp. 581–589, May. 2011.

[24] Pérez E, Ntaimo L, Ding Y. "Multi-component wind turbine modeling and simulation for wind farm operations and maintenance". *Simulation*, 91(4):360-82, Apr 2015.

[25] A. K. S. Jardine, D. Lin, and D. Banjevic, "A review on machinery diagnostics and prognostics implementing condition-based maintenance," *Mechanical Systems and Signal Processing*, vol.20, no.7, pp. 1483-1510, Oct. 2006.

[26] A. K. S. Jardine, T. Joseph, and D. Banjevic, "Optimizing condition-based maintenance decisions for equipment subject to vibration monitoring," *Journal of Quality in Maintenance Engineering*, vol. 5, no. 3, pp.192 – 202, 1995.

[27] Y. Zhu, E. Elsayed, H. Liao, and L. Chan, "Availability optimization of systems subject to competing risk," *European Journal of Operational Research*, vol. 202, no. 3, pp. 781-788, May. 2010.

[28] Y. Wang, H. Pham, "Multi-objective optimization of imperfect preventive maintenance policy for dependent competing risk system with hidden failure," *IEEE Transactions on Reliability*, vol. 60, no. 4, pp.770-781, Sept. 2011.

[29] Ko YM, Byon E. "Condition-based joint maintenance optimization for a large-scale system with homogeneous units". *IISE Transactions*. 49(5):493-504, May 2017.

[30] Abdul-Malak DT, Kharoufeh JP. "Optimally Replacing Multiple Systems in a Shared Environment". *Probability in the Engineering and Informational Sciences*:1-28, May 2017.

[31] Wang H, Wang W, Peng R. "A two-phase inspection model for a single component system with three-stage degradation". *Reliability Engineering & System Safety*,158:31-40, Feb 2017.

[32] Sung CK, Sheu SH, Hsu TS, Chen YC. "Extended optimal replacement policy for a two-unit system with failure rate interaction and external shocks". *International Journal of Systems Science*. 44(5):877-88, May 2013.

[33] Van Noortwijk JM. "A survey of the application of gamma processes in maintenance" *Reliability Engineering & System Safety*.,94(1):2-1, Jan 2009.

[34] Caballé NC, Castro IT, Pérez CJ, Lanza-Gutiérrez JM. "A condition-based maintenance of a dependent degradation-threshold-shock model in a system with multiple degradation processes". *Reliability Engineering & System Safety*.134:98-109, Feb 2015

[35] H. Peng, Q. Feng, and D. W. Coit, "Reliability and maintenance modeling for systems subject to multiple dependent competing failure processes," *IIE Transactions*, vol. 43, no. 1, pp.12-22, Apr. 2010.

[36] R. A. Waltz, J. L. Morales, J. Nocedal, and D. Orban, "An interior algorithm for nonlinear optimization that combines line search and trust region steps," *Mathematical Programming*, vol. 107, no. 3, pp. 391-408, Jul. 2006.
27